\documentclass[11pt,english]{article}

\usepackage[T1]{fontenc}
\usepackage[latin9]{inputenc}
\usepackage{babel}

\usepackage{booktabs}
\usepackage{url}
\usepackage{amsthm}
\usepackage{amsmath,amssymb}
\usepackage{graphicx}
\usepackage{setspace}
\usepackage[unicode=true]
 {hyperref}

\makeatletter

\providecommand{\tabularnewline}{\\}

\theoremstyle{plain}
\newtheorem{thm}{Theorem}
  \theoremstyle{definition}
  \newtheorem{defn}[thm]{Definition}
 \theoremstyle{definition}
  \newtheorem{example}[thm]{Example}
  \theoremstyle{plain}
  \newtheorem{lem}[thm]{Lemma}
  \theoremstyle{plain}
  \newtheorem{cor}[thm]{Corollary}

\usepackage{kbordermatrix}
\renewcommand{\kbldelim}{(}
\renewcommand{\kbrdelim}{)}
\usepackage[hang,small,bf]{caption}
\usepackage{framed}
\usepackage{xcolor}
\usepackage{fullpage}
\pdfoutput=1

\colorlet{shadecolor}{blue!25}   
\renewenvironment{leftbar}{%
  \MakeFramed {\advance\hsize-\width \FrameRestore}}%
{\endMakeFramed}

\theoremstyle{plain}
\newtheorem{mthm}{Theorem}

\newtheorem{mlemma}[mthm]{Lemma}
\newtheorem{mcor}[mthm]{Corollary}

\theoremstyle{definition}
\newtheorem{mdef}{Definition}
\newtheorem{mexample}{Example}
 
\renewenvironment{thm}
  {\begin{leftbar}\begin{mthm}}
  {\end{mthm}\end{leftbar}}
 
\renewenvironment{defn}
  {\begin{leftbar}\begin{mdef}}
  {\end{mdef}\end{leftbar}}

\renewenvironment{example}
  {\begin{leftbar}\begin{mexample}}
  {\end{mexample}\end{leftbar}}

\renewenvironment{cor}
{\begin{leftbar}\begin{mcor}}
{\end{mcor}\end{leftbar}}

\renewenvironment{lem}
  {\begin{leftbar}\begin{mlemma}}
  {\end{mlemma}\end{leftbar}}

\usepackage{tikz}
\usetikzlibrary{matrix,decorations.pathreplacing,calc}

\pgfkeys{tikz/mymatrixenv/.style={decoration=brace,every left delimiter/.style={xshift=0pt},every right delimiter/.style={xshift=0pt}}}
\pgfkeys{tikz/mymatrix/.style={matrix of math nodes,left delimiter=(,right delimiter={)},inner sep=2pt,column sep=1em,row sep=0.5em,nodes={inner sep=0pt}}}
\pgfkeys{tikz/mymatrixbrace/.style={decorate,thick}}

\newcommand\mymatrixbraceoffsetv{0.2em}

\newcommand*\mymatrixbracetop[4][m]{
    \draw[mymatrixbrace] ($(#1.north west)!(#1-1-#2.north west)!(#1.north east)+(0,\mymatrixbraceoffsetv)$)
        -- node[above=2pt] {#4} 
        ($(#1.north west)!(#1-1-#3.north east)!(#1.north east)+(0,\mymatrixbraceoffsetv)$);
}

\newcommand*\mymatrixnakedtop[4][m]{
    \draw[thin] ($(#1.north west)!(#1-1-#2.north west)!(#1.north east)+(0,\mymatrixbraceoffsetv)$)
        -- node[above=2pt] {#4} 
        ($(#1.north west)!(#1-1-#3.north east)!(#1.north east)+(0,\mymatrixbraceoffsetv)$);
}

\makeatother

\begin{document}

\title{The Kalmanson Complex\let\thefootnote\relax}

\author{Jonathan Terhorst%
\thanks{Department of Mathematics, San Francisco State University. Author
contact: \protect\href{mailto:terhorst@sfsu.edu}{terhorst@sfsu.edu}%
}}
\maketitle
\begin{abstract}
Let $X$ be a finite set of cardinality $n$. The Kalmanson complex
$\mathcal{K}_{n}$ is the simplicial complex whose vertices are non-trivial
$X$-splits, and whose facets are maximal circular split systems over
$X$. In this paper we examine $\mathcal{K}_{n}$ from three perspectives.
In addition to the $T$-theoretic description, we show that $\mathcal{K}_{n}$
has a geometric realization as the Kalmanson conditions on a finite
metric. A third description arises in terms of binary matrices which
possess the circular ones property. We prove the equivalence of these
three definitions. This leads to a simplified proof of the well-known
equivalence between Kalmanson and circular decomposable metrics, as
well as a partial description of the $f$-vector of \emph{$\mathcal{K}_{n}$.}

\end{abstract}
\newpage{}

\section{\label{sec:Introduction}Introduction}

A phylogenetic tree is a connected, acyclic graph which presents the
common evolutionary history of a group of species (taxa). A \emph{phylogenetic
network} generalizes this structure by allowing for the presence of
cycles. Phylogenetic networks have become a popular means of conveying
recombination, horizontal transfer and other reticulate events which
cannot be represented by a tree \cite{Arenas:2008uq,Huson:2006kx}. 

A particularly simple and well-known form of phylogenetic network
is the split network. As explained in \cite{Bryant:2004p1070}, split
networks are mathematically founded on $T$-theory (cf.\ \cite{Dress:1996p3275})
and express the so-called \emph{circular decomposition} of a finite
metric \cite{Bandelt:1992p3044,Christopher:1996p615}. The necessary
conditions for such a decomposition are a set of linear inequalities,
so that the space of circular decomposable metrics possesses polyhedral
structure.

In this paper we investigate that structure. Permuting these inequalities
produces a set of polyhedra whose union contains all circular decomposable
metrics. Since the polyhedra intersect along faces, the resulting
face lattice forms a simplicial complex. (We call this the Kalmanson
complex after \cite{Kalmanson:1975p495}, who derived the original
inequalities while studying certain tractable instances of the traveling
salesman problem.) 

Abstracting away from the underlying geometry, we show how combinatorially
isomorphic objects can be derived in terms of either circular split
systems, or binary matrices which possess the consecutive ones property
\cite{Booth:1976p5201}. Our main result is to prove the equivalence
of these structures by exhibiting order-preserving bijections between
their face lattices (Theorems \ref{thm:K-fan iso} and \ref{thm:circular iff Circ1R}).
This in turn leads to a new proof of the equivalence of Kalmanson
and circular decomposable metrics (Corollary \ref{cor:circdecompiffkalmanson}).
This is a known result \cite{Chepoi:1998p534,Christopher:1996p615},
but our proof has the benefit of being extremely simple, relying only
on basic concepts from polyhedral geometry.

We then use these findings to study the $f$-vector of Kalmanson complex.
We relate the problem of enumerating its faces to a counting problem
on certain classes of binary matrices, and exploit a structure theorem
of \cite{Tucker:1972p5021} to obtain a new result on the number of
triangles contained (Theorem \ref{thm:triangles}). Even in the simplest
non-trivial case, this counting problem is seen to possess considerable
complexity, and we leave a more general method of counting the faces
of the Kalmanson complex as an interesting open problem.

The paper is organized as follows. Section \ref{sec:Definitions}
begins with some preliminaries from $T$-theory which enable us to
define the complex abstractly in terms of split systems. In Section
\ref{sec:Geometry-of}, we review the Kalmanson conditions. These
are a set of inequality restrictions on a finite metric which, when
satisfied, allow the traveling salesman problem to be solved in constant
time. We show that these inequalities give a geometric realization
of the Kalmanson complex. In Section \ref{sec:-and-the}, we study
the consecutive ones property for binary matrices. This property is
shown to be equivalent to the circularity property for split systems
discussed above, giving us a third description of the complex in terms
of equivalence classes of binary matrices. In Section \ref{sec:-vector},
we use these three ways of viewing the Kalmanson complex to enumerate
some of its faces, thus giving a partial characterization of its $f$-vector.
Finally, in Section \ref{sec:Conclusion} we offer some concluding
remarks.

\paragraph*{Acknowledgments.}

This material is adapted from my master's thesis at SFSU. I am indebted
to my adviser, Serkan Ho\c{s}ten, for his patience, encouragement
and support. I would also like to thank Bernd Sturmfels, David Bryant,
Lior Pachter, Alex Engström and Raman Sanyal for their helpful suggestions
and comments. Of course, all errors are my own.

\section{\label{sec:Definitions}Definitions}

We begin with some basic concepts from $T$-theory. For a full introduction,
see \cite{Bandelt:1992p3044,Dress:1996p3275}. 

Throughout the paper, $X$ is a finite set of cardinality $n\geq4$.
An \emph{$X$-split }is a bipartition of $X$; that is, $S=\{A,B\}$
is an $X$-split if $A\cap B=\emptyset$ and $A\cup B=X$. (When the
meaning is obvious, we will simply call $S$ a split.) $A$ and $B$
are called the \emph{blocks} of $S$, and the \emph{size }of $S$
is defined as $\operatorname{size}(S):=\min\left\{ \left|A\right|,\left|B\right|\right\} $.
$S$ is \emph{non-trivial} if $\operatorname{size}(S)>1$ and \emph{minimal}
if $\operatorname{size}\left(S\right)=2$.

Let $\mathcal{S}(X)$ be the set of non-trivial $X$-splits. A \emph{split
system }$\mathcal{S}\subset\mathcal{S}(X)$ is a set of splits. \cite{Bandelt:1992p3044}
introduced the concept of a \emph{circular split system.}%
\footnote{Circular split systems are sometimes referred to in the literature
as \emph{cyclic split systems}.%
}
\begin{defn}
\label{def:circular split system}A split system $\mathcal{S}$ is
\emph{circular} if there is a permutation $\sigma\in S_{n}$ such
that for each split $S=\{A,B\}\in\mathcal{S}$ there exists $i,j\in[n]$
such that \[
S=\Bigg\{\left\{ x_{\sigma\left(\overline{i}\right)},x_{\sigma\left(\overline{i+1}\right)},\dots,x_{\sigma\left(\overline{j-1}\right)},x_{\sigma\left(\overline{j}\right)}\right\} ,\left\{ x_{\sigma\left(\overline{j}\right)},x_{\sigma\left(\overline{j+1}\right)},\dots,x_{\sigma\left(\overline{i-1}\right)},x_{\sigma\left(\overline{i}\right)}\right\} \Bigg\}\]
where $\overline{i}$ denotes $i\pmod n$.
\end{defn}
\begin{figure}
\noindent \begin{centering}
\includegraphics[bb=0bp 0bp 335bp 191bp,scale=0.7]{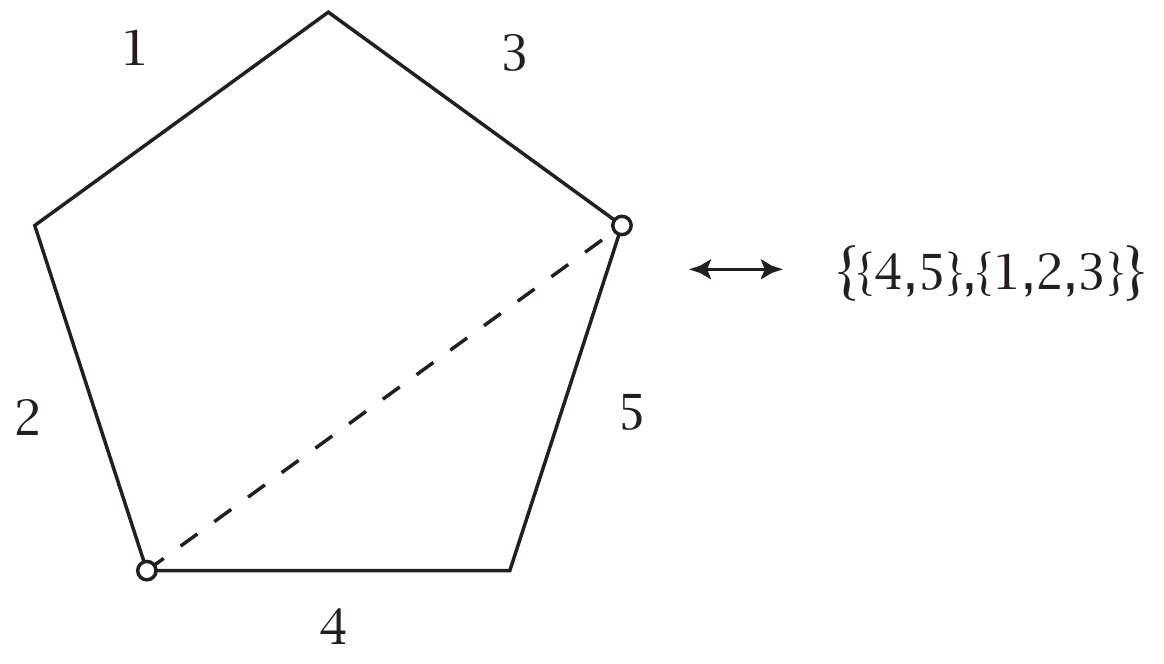}
\par\end{centering}

\caption{\label{fig:circular-split}A circular split.}
\end{figure}

Circular split systems have a simple geometric interpretation: they
are obtained by labeling the edges of a regular $n$-gon, and connecting
its edges with diagonals to form splits (Figure \ref{fig:circular-split}).
From this we see that a circular split system contains at most $\binom{n}{2}$
distinct splits. A partial converse also holds: a weakly compatible
(cf. Section \ref{sec:-vector}) split system containing $\binom{n}{2}$
splits is circular \cite{Bandelt:1992p3044}.

From Definition \ref{def:circular split system} we see the set of
circular split systems is closed under the operations of taking subsets
(any subset of a circular split system is circular) and forming intersections.
Hence, it is a simplicial complex. This complex is our main object
of study.
\begin{defn}
The \emph{Kalmanson complex} is the simplicial complex whose vertices
are $X$-splits, and whose facets are maximal circular split systems.
\end{defn}
Clearly this complex is unique up to the cardinality of $X$. Henceforth
we write $\mathcal{K}_{n}$ to denote the Kalmanson complex over a
base set of cardinality $n$.

\section{\label{sec:Geometry-of}Geometry of $\mathcal{K}_{n}$}

Throughout this section, we let $D=\left(d_{ij}\right)_{i,j\in[n]}$
be a symmetric, non-negative matrix with zeros along the diagonal.
We refer to matrices possessing this property as \emph{distance matrices.}
\begin{defn}
Let $D=\left(d_{ij}\right)$ be a distance matrix and let $S_{n}$
denote the symmetric group on $n$ letters. The \emph{traveling salesman
problem} (TSP) \emph{over $D$} is \[
\min_{\sigma\in S_{n}}\left(\sum_{i=1}^{n-1}d_{\sigma(i)\sigma(i+1)}+d_{\sigma(n)}{}_{\sigma(1)}\right)\]

\end{defn}
For general $D$, it is well-known that the TSP is NP-hard. However,
some special cases have lower complexity. In particular, \cite{Kalmanson:1975p495}
showed that if $D$ satisfies a certain set of linear inequalities,
then the TSP over $D$ possesses a trivial solution.
\begin{thm}[\cite{Kalmanson:1975p495}]
 Let $D$ be a distance matrix. If \begin{equation}
\max(d_{ij}+d_{kl},d_{il}+d_{jk})\leq d_{ik}+d_{jl}\ \mathcal{\text{for all}}\ 1\leq i<j<k<l\leq n\label{eq:kalmanson condition}\end{equation}
then the identity permutation solves the TSP over $D$.
\end{thm}
The inequalities \eqref{eq:kalmanson condition} are referred to as
the \emph{Kalmanson conditions}, and a matrix which satisfies them
is a \emph{Kalmanson matrix }(or simply Kalmanson.) 

It may be that $D$ does not satisfy \eqref{eq:kalmanson condition},
but that some permutation of the rows and columns of $D$ does. In
this case we say that $D$ is a \emph{permuted Kalmanson matrix}.
Since permuting $D$ amounts to simply relabeling the underlying distance
or cost data, this operation preserves the structure of the problem.
\cite{Christopher:1996p615} give an $O(n^{2})$ recognition algorithm
for permuted Kalmanson matrices, so we say the TSP is polynomial time-solvable
for this class.

Geometrically, \eqref{eq:kalmanson condition} comprises a finite
intersection of closed half-spaces: a polyhedron. Given a polyhedron
$P\subset\mathbb{R}^{k}$ and a hyperplane $H\subset\mathbb{R}^{k}$,
we say $H$ \emph{supports }$P$ if $H\cap P\neq\emptyset$ and $P$
is completely contained in one of the closed half-spaces defined by
$H$. $F\subset P$ is a \emph{face} of $P$ if $F=P\cap H$ for some
supporting hyperplane $H$ of $P$. The \emph{face lattice }of $P$
is the poset of faces of $P$ ordered by set inclusion.

Recall that a set of polyhedra which intersect along faces is called
a \emph{polyhedral fan}. Permuting the indices in \eqref{eq:kalmanson condition}
generates a polyhedral fan which we denote $\mathcal{P}_{n}$.
\begin{example}
For $n=4$, $\mathcal{P}_{n}$ is the union of three polyhedra obtained
by permuting the indices in \eqref{eq:kalmanson condition}: $\mathcal{P}_{4}:=\left(d_{ij}\right)_{i,j\in[4]}$
such that \[
\left\{ \begin{array}{cc}
d_{12}+d_{34} & \leq d_{14}+d_{23}\\
d_{13}+d_{24} & \leq d_{14}+d_{23}\end{array}\right\} \text{ or }\left\{ \begin{array}{cc}
d_{13}+d_{24} & \leq d_{12}+d_{34}\\
d_{14}+d_{23} & \leq d_{12}+d_{34}\end{array}\right\} \text{ or }\left\{ \begin{array}{cc}
d_{14}+d_{23} & \leq d_{13}+d_{24}\\
d_{12}+d_{34} & \leq d_{13}+d_{24}\end{array}\right\} \]
Collectively, these define the region of $\mathbb{R}^{\binom{4}{2}}$
containing all $4\times4$ permuted Kalmanson matrices.
\end{example}

\subsection{Equivalence of $\mathcal{P}_{n}$ and $\mathcal{K}_{n}$}

The main claim of this section is that $\mathcal{P}_{n}$ is a geometric
realization of $\mathcal{K}_{n}$ in the sense that they are combinatorially
equivalent.
\begin{thm}
\label{thm:K-fan iso}The face lattices of $\mathcal{K}_{n}$ and
$\mathcal{P}_{n}$ are isomorphic as posets.
\end{thm}
The remainder of this subsection is devoted to proving the theorem
by finding an inclusion-preserving bijection between the faces of
these two sets. 

In \cite{demidenko} it is shown that the polyhedron defined by \eqref{eq:kalmanson condition}
decomposes into an $n$-dimensional lineality space and a pointed
cone of dimension $\binom{n}{2}-n$. We are interested in the structure
of the latter since it encapsulates the combinatorial data embodied
by the polyhedron. The authors give an explicit description of the
extreme rays of this cone.
\begin{example}
\label{exa:rays}For $n=5$, the rays of the standard Kalmanson polyhedron
are \begin{align*}
V^{(2)} & =\begin{pmatrix}\begin{array}{rr|rrr}
0 & 0 & 1 & 1 & 1\\
0 & 0 & 1 & 1 & 1\\
\hline 1 & 1 & 0 & 0 & 0\\
1 & 1 & 0 & 0 & 0\\
1 & 1 & 0 & 0 & 0\end{array}\end{pmatrix} & V^{(3)} & =\begin{pmatrix}\begin{array}{rrr|rr}
0 & 0 & 0 & 1 & 1\\
0 & 0 & 0 & 1 & 1\\
0 & 0 & 0 & 1 & 1\\
\hline 1 & 1 & 1 & 0 & 0\\
1 & 1 & 1 & 0 & 0\end{array}\end{pmatrix}\\
V^{(1,3)} & =\begin{pmatrix}\begin{array}{r|rr|rr}
0 & 1 & 1 & 0 & 0\\
\hline 1 & 0 & 0 & 1 & 1\\
1 & 0 & 0 & 1 & 1\\
\hline 0 & 1 & 1 & 0 & 0\\
0 & 1 & 1 & 0 & 0\end{array}\end{pmatrix} & V^{(1,4)} & =\begin{pmatrix}\begin{array}{r|rrr|r}
0 & 1 & 1 & 1 & 0\\
\hline 1 & 0 & 0 & 0 & 1\\
1 & 0 & 0 & 0 & 1\\
1 & 0 & 0 & 0 & 1\\
\hline 0 & 1 & 1 & 1 & 0\end{array}\end{pmatrix} & V^{(2,4)} & =\left(\begin{array}{rr|rr|r}
0 & 0 & 1 & 1 & 0\\
0 & 0 & 1 & 1 & 0\\
\hline 1 & 1 & 0 & 0 & 1\\
1 & 1 & 0 & 0 & 1\\
\hline 0 & 0 & 1 & 1 & 0\end{array}\right)\end{align*}

\end{example}
We see that the $V^{(i)}$ and $V^{(i,j)}$ have a structure which
is the result of arranging square blocks of zeros along the diagonal.
It turns out that these matrices, along with their permutations, encode
the non-trivial $X$-splits.
\begin{defn}
Let $S=\left\{ A,B\right\} $ be an $X$-split. A \emph{split metric}
$\delta_{S}:X\times X\to\mathbb{R}$ is a function such that \[
\delta_{S}(x,y)=\begin{cases}
0, & \left\{ x,y\right\} \subset A\text{ or }\left\{ x,y\right\} \subset B\\
1, & \text{otherwise}\end{cases}\]

\end{defn}
Split metrics are unique.
\begin{lem}
Let $S_{1},S_{2}\in\mathcal{S}(X)$. If $\delta_{S_{1}}=\delta_{S_{2}},$
then $S_{1}=S_{2}$.
\begin{proof}
For a split $S$ define $\gamma_{S}(i):=\left\{ x\in X:\delta_{S}(1,x)=i\right\} $.
We have $S=\left\{ A,B\right\} =\left\{ \gamma_{S}(0),\gamma_{S}(1)\right\} $.
Hence $S_{1}=\left\{ \gamma_{S_{1}}(0),\gamma_{S_{1}}\left(1\right)\right\} =\left\{ \gamma_{S_{2}}(0),\gamma_{S_{2}}\left(1\right)\right\} =S_{2}$.
\end{proof}
\end{lem}
Returning to the example, let $\Delta_{S}$ be the (symmetric, $n\times n$)
matrix associated to $\delta_{S}$. The matrices in Example \ref{exa:rays}
are obtained from split metrics: \begin{align*}
V^{(2)} & =\Delta_{12|345} & V^{(3)} & =\Delta_{123|45}\\
V^{(1,3)} & =\Delta_{145|23} & V^{(1,4)} & =\Delta_{14|235} & V^{(2,4)} & =\Delta_{125|34}\end{align*}

We now formalize this idea.
\begin{thm}[\cite{demidenko}]
 \label{thm:V matrices}The space of Kalmanson matrices consists
of an $n$-dimensional lineality space and an $n(n-3)/2$-dimensional
pointed cone. The lineality space is spanned by the matrices $E^{(i)}=\left(e_{pq}^{(i)}\right),1\leq i\leq n$
where \begin{equation}
e_{pq}^{(i)}=\begin{cases}
1, & p=i\text{ xor }q=i\\
0, & \text{otherwise}\end{cases}\label{eq:eij}\end{equation}
The pointed cone is ruled by the symmetric matrices $V^{(i)}=\left(v_{pq}^{(i)}\right)$,
\textup{$2\leq i\le n-2$} and $V^{(i,j)}=\left(v_{pq}^{(i,j)}\right)$,
$1\leq i\leq n-3$, $i+2\leq j\leq n-1$, where \begin{align}
v_{pq}^{(i)} & :=\begin{cases}
1, & 1\leq p\le i<q\leq n\\
0, & \text{otherwise }\end{cases}\label{eq:vi}\\
v_{pq}^{(i,j)} & :=\begin{cases}
1, & 1\leq p\leq i<q\leq j\text{ or }i<p\leq j<q\leq n\\
0, & \text{otherwise}\end{cases}\label{eq:vij}\end{align}

\end{thm}
Now let the symmetric group $S_{n}$ act on the set of $n\times n$
matrices by symmetric permutation of rows and columns: $\sigma\cdot M=\left(m_{\sigma(i),\sigma(j)}\right)$
for all $\sigma\in S_{n}$ and $M=\left(m_{ij}\right)$. The following
lemma is immediate from \eqref{eq:vi} and \eqref{eq:vij}.
\begin{lem}
\label{lem:splitmetrics}$V^{(i)}=\Delta_{1\cdots i|i+1\cdots n}$
and $V^{(i,j)}=\Delta_{i+1\cdots j|j+1\cdots i}$. Additionally, symmetrically
permuting $V^{(i)},V^{(i,j)}$ is equivalent to applying the same
permutation to the underlying split: $\sigma\cdot V^{(i)}=\Delta_{\sigma(1)\cdots\sigma(i)|\sigma(i+1)\cdots\sigma(n)}$
(and similarly for $V^{(i,j)}$).
\end{lem}
Define \begin{align*}
\mathcal{V} & :=\left\{ V^{(i)}:2\leq i\leq n-2\right\} \cup\left\{ V^{(i,j)}:1\leq i\leq n-3,i+2\leq j\leq n-1\right\} \\
\mathcal{R} & :=\left\{ \sigma\cdot V:\sigma\in S_{n},V\in\mathcal{V}\right\} \end{align*}
Thus $\mathcal{R}$ is the set of vertices (rays) of $\mathcal{P}_{n}$.
Finally, let $T(V)$ be the map which takes a matrix in $\mathcal{R}$
to its corresponding split,\begin{align*}
T:\mathcal{R} & \to\mathcal{S}(X)\\
\sigma\cdot V^{(i)} & \mapsto\Big\{\left\{ \sigma(1),\dots,\sigma(i)\right\} ,\left\{ \sigma(i+1),\dots,\sigma(n)\right\} \Big\}\end{align*}

(Note that $T$ is well-defined since for each $i,j\in[n]$ there
exists a $\sigma\in S_{n}$ such that $V^{(i,j)}=\sigma\cdot V^{(j-i)}$.)
\begin{lem}
$T:\mathcal{R}\to\mathcal{S}(X)$ is a bijection.
\begin{proof}
Injectivity follows from the uniqueness of split metrics. For surjectivity,
let $T=\left\{ A,B\right\} \in\mathcal{S}(X)$ be a split. Let $\sigma\in S_{n}$
be a permutation such that $\sigma^{-1}\cdot A=\left\{ 1,\dots,|A|\right\} $.
Then $T\left(\sigma\cdot V^{(|A|)}\right)=S$.
\end{proof}
\end{lem}
It remains to show that $T$ is order-preserving: $T(U)\subset T(V)\iff U\subset V$.
This follows from the fact that $T$ maps faces to faces.
\begin{lem}
$\left\{ T(M_{1}),T(M_{2}),\dots,T(M_{k})\right\} $ is a face of
$\mathcal{K}_{n}$ if and only if $\left\{ M_{1},M_{2},\dots,M_{k}\right\} $
is a face of $\mathcal{P}_{n}$.
\begin{proof}
If $\left\{ T(M_{1}),T(M_{2}),\dots,T(M_{k})\right\} $ is a face
of $\mathcal{K}_{n}$ then it is circular respect to the ordering
$\left(\sigma(1),\dots,\sigma(n)\right)$ for some $\sigma\in S_{n}$.
Hence $\left\{ M_{1},M_{2},\dots,M_{k}\right\} \subseteq\sigma\cdot\mathcal{V}$
is a face of $\mathcal{P}_{n}$.

Conversely, since \begin{align*}
T(V^{(i)}) & =\Big\{\left\{ 1,2,\dots,i\right\} ,\left\{ i+1,\dots,n\right\} \Big\}\\
T(V^{(i,j)}) & =\Big\{\left\{ 1,2,\dots,i,j+1,\dots,n\right\} ,\left\{ i+1,\dots,j\right\} \Big\}\end{align*}
the set $S:=T(\mathcal{V})$ is a maximal circular split system with
the ordering $(1,2,\dots,n)$. Therefore the claim is true when the
$M_{i}\in\mathcal{V}$. 

Now if $M_{1},\dots,M_{k}$ is an arbitrary face of $\mathcal{P}_{n}$,
then there is a $\sigma\in S_{n}$ such that for each $1\leq i\le k$
there exists $V_{i}\in\mathcal{V}$ with $M_{i}=\sigma\cdot V_{i}$.
Then \begin{align*}
\left\{ T(M_{1}),\dots,T(M_{k})\right\}  & =\left\{ T(\sigma\cdot V_{1}),\dots,T(\sigma\cdot V_{k})\right\} \\
 & =\left\{ \sigma\cdot T(V_{1}),\dots,\sigma\cdot T(V_{k})\right\} \\
 & \subseteq\sigma\cdot S\end{align*}
 is a face of $\mathcal{K}_{n}$.
\end{proof}
\end{lem}
This concludes the proof of Theorem \ref{thm:K-fan iso}.

\subsection{Circular Decomposability}

A metric $\delta$ is \emph{circular decomposable} if it can be written
as the positively-weighted sum of circular split metrics, i.e.\[
\delta=\sum_{S\in\mathcal{C}}\alpha_{S}\delta_{S}\]
for some circular split system $\mathcal{C}$ and weights $\alpha_{S}>0$.
It has been shown in \cite{Chepoi:1998p534} and \cite{Christopher:1996p615}
that $\delta$ is circular decomposable if and only if it satisfies
the Kalmanson conditions. Both proofs are non-trivial; \cite{Chepoi:1998p534}
relies on a Crofton-type formula for computing distances in metric
spaces, while \cite{Christopher:1996p615} uses a number of results
from the theory of metrics over a finite set \cite{Bandelt:1992p3044}.

The polyhedral characterization of the Kalmanson cone given in \cite{demidenko},
coupled with the observations of the preceding section, enable us
to establish this equivalence in a new and straightforward way. Indeed,
$\delta$ satisfies the Kalmanson conditions iff it is in a permuted
Kalmanson cone. By Theorem \ref{thm:V matrices} this occurs iff $\delta$
is a linear combination of the permuted matrices $E^{(i)},V^{(i)},V^{(i,j)}$
for some $\sigma\in S_{n}$:\begin{equation}
\delta=\sigma\cdot\left(\sum_{i=1}^{n}\alpha_{i}E^{(i)}+\sum_{i=2}^{n-2}\beta_{i}V^{(i)}+\sum_{i=1}^{n-3}\sum_{j=i+2}^{n-1}\gamma_{ij}V^{(i,j)}\right)\text{ for }\alpha_{i}\in\mathbb{R}\text{ and }\beta_{i},\gamma_{ij}>0\label{eq:delta}\end{equation}
By Lemma \ref{lem:splitmetrics}, $V^{(i)}$ and $V^{(i,j)}$ are
split metrics, and it is easily seen that the $E^{(i)}$ are split
metrics corresponding to trivial splits.

Now, equation \eqref{eq:delta} is not necessarily a circular decomposition
since the $\alpha_{i}$ can be negative. However, assuming $\delta$
obeys the triangle inequality (which, recall, is not implied by the
Kalmanson conditions), the $\alpha_{i}$ are seen to be\emph{ }non-negative.
\begin{lem}
Let $\delta$ be written as in \eqref{eq:delta}. Then \[
\delta(i,i+1)+\delta(i+1,i+2)-\delta(i,i+2)=2\alpha_{i+1}\]
where $i,i+1,i+2$ are modulo $n$.
\begin{proof}
Put $\delta=\delta_{\alpha}+\delta_{\beta}+\delta_{\gamma}$, where
\begin{align*}
\delta_{\alpha} & =\sum_{i=1}^{n}\alpha_{i}E^{(i)} & \delta_{\beta} & =\sum_{i=2}^{n-2}\beta_{i}V^{(i)} & \delta_{\gamma} & =\sum_{i=1}^{n-3}\sum_{j=i+2}^{n-1}\gamma_{ij}V^{(i,j)}\end{align*}
From equations \eqref{eq:eij}--\eqref{eq:vij} we have \begin{align}
\delta_{\alpha}(i,j) & =\alpha_{i}+\alpha_{j}\\
\delta_{\beta}(i,j) & =\sum_{\substack{2\leq s\leq n-2\\
i\le s<j}
}\beta_{s}\\
\delta_{\gamma}(i,j) & =\sum_{\substack{i\leq s<j\leq t\\
1\leq s\leq n-3\\
s+2\leq t\leq n-1}
}\gamma_{st}\;+\sum_{\substack{s<i\le t<j\\
1\leq s\leq n-3\\
s+2\leq t\leq n-1}
}\gamma_{st}\label{eq:delta_elements}\end{align}
Hence $\delta_{\alpha}(i,i+1)+\delta_{\alpha}(i+1,i+2)-\delta_{\alpha}(i,i+2)=2\alpha_{i+1}$. 

For $\delta_{\beta}$ we have \begin{align*}
\delta_{\beta}(i,i+1) & =\begin{cases}
0, & i=1,n-1\\
\beta_{i}, & \text{otherwise}\end{cases}\\
\delta_{\beta}(i,i+2) & =\begin{cases}
\beta_{i}, & i=1,n-2\\
\beta_{i}+\beta_{i+1}, & \text{otherwise}\end{cases}\\
\delta_{\beta}(1,n)=\delta_{\beta}(2,n) & =\beta_{2}+\cdots+\beta_{n-2}\end{align*}
We see that $\delta_{\beta}(i,i+1)+\delta_{\beta}(i+1,i+2)=\delta_{\beta}(i,i+2)$. 

Finally, for $\delta_{\gamma}$ we further decompose it as $\delta_{\gamma}=\delta_{\gamma_{1}}+\delta_{\gamma_{2}}$
according to the two summands in \eqref{eq:delta_elements}. Repeating
the same procedure yields \begin{align*}
\delta_{\gamma_{1}}(i,i+1) & =\sum_{t=i+2}^{n-1}\gamma_{i,t}\\
\delta_{\gamma_{1}}(i,i+2) & =\sum_{a=0}^{1}\:\sum_{t=i+2+a}^{n-1}\gamma_{i+a,t}\\
\delta_{\gamma_{1}}(1,n)=\delta_{\gamma_{2}}(2,n) & =0\\
\delta_{\gamma_{2}}(i,i+1)=\delta_{\gamma_{2}}(i,i+2)=\delta_{\gamma_{2}}(1,n) & =0\\
\delta_{\gamma_{2}}(2,n) & =\sum_{t=3}^{n-1}\gamma_{1,t}\end{align*}
After some algebraic manipulations we again obtain $\delta_{\gamma}(i,i+1)+\delta_{\gamma}(i+1,i+2)=\delta_{\gamma}(i,i+2)$.
\end{proof}
\end{lem}
\begin{cor}
\label{cor:circdecompiffkalmanson}Let $\delta:X\times X\to\mathbb{R}$
be a metric over the finite set $X$. Then $\delta$ is circular decomposable
if and only if it satisfies the Kalmanson conditions.
\end{cor}

\section{\label{sec:-and-the}$\mathcal{K}_{n}$ and the consecutive ones
property}

Thus far we have defined $\mathcal{K}_{n}$ as a split-theoretic simplicial
complex and also geometrically in terms of permutations of the Kalmanson
conditions. In this section we present a third description of the
Kalmanson complex as a set of (equivalence classes of) binary matrices
possessing a certain structure. Again, we will show that this formulation
is entirely equivalent to the preceding two. Throughout this section,
$ $$M$ is taken to be an $m\times n$ binary matrix (entries are
zero or one.) 
\begin{defn}
\label{def:c1r}$M$ is said to possess the \emph{consecutive ones
property for rows} (C1R) if its columns may be permuted such that
the ones in each row occur in blocks. $M$ possesses the \emph{circular
ones property for rows} (Circ1R) if its columns may be permuted such
that either the ones or the zeros (or both) in each row occur in a
block.
\end{defn}
Intuitively, a Circ1R matrix has the property that for each of its
rows, the ones occur in a block when it is {}``wrapped around''
a cylinder.
\begin{example}
Consider the matrices\[
\underset{(1)}{\begin{pmatrix}1 & 1 & 0 & 0 & 0\\
0 & 0 & 1 & 1 & 0\\
1 & 0 & 0 & 0 & 0\\
0 & 1 & 1 & 1 & 0\end{pmatrix}}\qquad\underset{(2)}{\begin{pmatrix}0 & 0 & 1 & 1 & 0\\
1 & 1 & 1 & 0 & 0\\
0 & 1 & 0 & 0 & 1\\
1 & 1 & 1 & 0 & 1\end{pmatrix}}\qquad\underset{(3)}{\begin{pmatrix}1 & 0 & 0 & 1 & 1\\
0 & 1 & 1 & 0 & 0\\
1 & 1 & 0 & 1 & 1\\
0 & 0 & 0 & 0 & 1\end{pmatrix}}\]
$(1)$ and $(2)$ are C1R, and $(3)$ is Circ1R. To verify that $(2)$
is C1R, we apply the permutation $(1\,3\,4\,5)\in S_{5}$ to its columns:

\[
\begin{pmatrix}0 & 0 & 1 & 1 & 0\\
1 & 1 & 1 & 0 & 0\\
0 & 1 & 0 & 0 & 1\\
1 & 1 & 1 & 0 & 1\end{pmatrix}\overset{(1\,3\,4\,5)}{\longrightarrow}\begin{pmatrix}0 & 0 & 0 & 1 & 1\\
0 & 1 & 1 & 1 & 0\\
1 & 1 & 0 & 0 & 0\\
1 & 1 & 1 & 1 & 0\end{pmatrix}\]

\end{example}
If $M$ is C1R or Circ1R, then the matrix obtained by replacing any
number of rows of $M$ by their binary complement will be Circ1R.
This provides justification for the following theorem.
\begin{thm}[\cite{Tucker:1970p4803}]
\label{lem:circ1r <-> c1r}Let $M$ be a binary matrix, and let $M'$
be the matrix obtained by complementing each row in $M$ which has
a one in the first column. Then $M$ is C1R if and only if $M'$ is
Circ1R.
\end{thm}
A circular split system and a Circ1R binary matrix are, in a sense,
identical. To see this, let $m$ be fixed and consider the set of
all split systems over $X$ which contain $m$ splits: $\mathcal{S}_{m}(X)=\left\{ \mathcal{S}\subset\mathcal{S}(X):|\mathcal{S}|=m\right\} $.
Also, let $\mathcal{M}_{m\times n}^{0}\left(\left\{ 0,1\right\} \right)$
be the set of $m\times n$ binary matrices who first column contains
all zeros, and let the symmetric group $S_{m}$ act on it by permutation
of rows. 

Finally, let $\mathcal{Q}_{m}=\mathcal{M}_{m\times n}^{0}\left(\left\{ 0,1\right\} \right)/\mathord{\sim}$
be the set of equivalence classes under the relation {}``$M_{1}\sim M_{2}\iff M_{1}=\sigma\cdot M_{2}$
for some $\sigma\in S_{m}$''. (Note that, as row permutations do
nothing to affect the C1R/Circ1R properties, it makes sense to say
that a class $[M]\in\mathcal{Q}$ possesses one or both.)

Now define a map $F:\mathcal{S}_{m}(X)\to\mathcal{Q}_{m}$ which sends
a system of $m$ splits to the class of the binary matrix obtained
by converting the splits to a binary vector and stacking them. Formally,\begin{align*}
F:\mathcal{S}_{m}(X) & \to\mathcal{Q}_{m}\\
\Big\{\left\{ A_{1},B_{1}\right\} ,\dots,\left\{ A_{m},B_{m}\right\} \Big\} & \mapsto\left[\left(w_{ij}\right)_{i\in[m],j\in[n]}\right]\\
w_{ij} & :=\begin{cases}
1, & j\in A_{i}\text{ and }1\notin A_{i}\\
1, & j\notin A_{i}\text{ and }1\in A_{i}\\
0, & \text{otherwise}\end{cases}\end{align*}

\begin{example}
Let $n=5$ and $\mathcal{S}\in\mathcal{S}_{3}(X)$ be the split system
\[
\mathcal{S}=\biggl\{\Bigl\{\left\{ 1,2\right\} ,\left\{ 3,4,5\right\} \Bigr\},\Bigl\{\left\{ 1,3,5\right\} ,\left\{ 2,4\right\} \Bigr\},\Bigl\{\left\{ 1,4\right\} ,\left\{ 2,3,5\right\} \Bigr\}\biggr\}\]
Then \[
F(S)=\left[\begin{pmatrix}0 & 0 & 1 & 1 & 1\\
0 & 1 & 0 & 1 & 0\\
0 & 1 & 1 & 0 & 1\end{pmatrix}\right]=\left[\begin{pmatrix}0 & 1 & 1 & 0 & 1\\
0 & 1 & 0 & 1 & 0\\
0 & 0 & 1 & 1 & 1\end{pmatrix}\right]\]
\end{example}
\begin{lem}
\label{lem:sm bijection}$F:\mathcal{S}_{m}(X)\to\mathcal{Q}_{m}$
is a bijection.
\begin{proof}
Let $[M]\in\mathcal{Q}_{m}$ be given. Simply convert each row of
$M\in\mathcal{M}_{m\times n}^{0}\left(\left\{ 0,1\right\} \right)$
to an $X$-split in the obvious way. The resulting split system \emph{$\mathcal{S}$}
gives $F(\mathcal{S})=[M]$, so $F$ is onto.

Now suppose $F(\mathcal{S}_{1})=F(\mathcal{S}_{2})$ for two split
systems $\mathcal{S}_{1},\mathcal{S}_{2}$. Then $\mathcal{S}_{1}$
and $\mathcal{S}_{2}$ represent the same splits up to ordering. But
then $\mathcal{S}_{1}=\mathcal{S}_{2}$ so $F$ is one-to-one.
\end{proof}
\end{lem}
Having established the bijection, it is easy to see that Circ1R and
circularity are analogous properties for $\mathcal{Q}$ and $\mathcal{S}_{m}(X)$,
respectively.
\begin{thm}
\label{thm:circular iff Circ1R}Let $\mathcal{S}\in\mathcal{S}_{m}(X)$
be an arbitrary split system. Then $\mathcal{S}$ is circular iff
$F(\mathcal{S})$ is Circ1R.
\begin{proof}
A split system $\mathcal{S}$ is circular iff there is a $\sigma\in S_{n}$
such that each split $S\in\mathcal{S}$ is of the form\[
S=\Bigg\{\left\{ \sigma(\overline{i}),\sigma(\overline{i+1}),\dots,\sigma(\overline{j})\right\} ,\left\{ \sigma(\overline{j+1}),\dots,\sigma(\overline{i-1})\right\} \Bigg\}\]
where $\overline{i}$ denotes $i\pmod n$. This occurs iff applying
$\sigma$ to the columns of $F(S)$ yields a class of matrices $[M]\in\mathcal{Q}_{m}$
whose ones appear consecutively. Since the first column of $M$ is
the zero vector, $M$ is C1P iff it is Circ1R by Theorem \ref{lem:circ1r <-> c1r}.
\end{proof}
\end{thm}
\begin{cor}
\label{cor:circular iff c1r}$\mathcal{S}\in\mathcal{S}_{m}(X)$ is
circular iff $F(\mathcal{S})$ is C1R.
\end{cor}
The preceding theorem enables us to furnish another description of
$\mathcal{K}_{n}$: it is the poset of all Circ1R binary matrices
(up to row permutation) which possess an initial column of zeros,
and at least two ones and two zeros in each row, ordered by inclusion
of the set of row vectors corresponding to each matrix.

\section{\label{sec:-vector}$f$-vector}

In this section we will harness the three descriptions of $\mathcal{K}_{n}$
to study its combinatorial structure in greater detail.
\begin{defn}
Let $\Delta$ be a simplicial complex of dimension $d-1$, and let
$f_{i}$ denote the number of $i$-dimensional faces of $\Delta$.
The \emph{$f-$vector of $\Delta$ }is the vector\emph{ $f=\left(f_{0},f_{1},\dots,f_{d-1}\right)$.}
\end{defn}
Thus, $f_{0}$ counts the vertices of $\Delta$ and $f_{d-1}$ counts
the facets. By geometric analogy, $f_{1}$, $f_{2}$ and $f_{d-2}$
are called the \emph{edges}, \emph{triangles}, and \emph{ridges} of
$\Delta$, respectively. 

\begin{table}
\noindent \begin{centering}
\begin{tabular}{cc}
\toprule 
\addlinespace
$n$ & $f$-vector\tabularnewline\addlinespace
\midrule
\addlinespace
\addlinespace
4 & $\left\langle 3,3\right\rangle $\tabularnewline\addlinespace
\addlinespace
\addlinespace
5 & $\left\langle 10,45,90,60,12\right\rangle $\tabularnewline\addlinespace
\addlinespace
\addlinespace
$6$ & $\left\langle 25,300,1755,4725,6390,4860,2160,540,60\right\rangle $\tabularnewline\addlinespace
\addlinespace
\addlinespace
7 & $\left\langle 56,1540,19950,121485,\dots,5040,360\right\rangle $\tabularnewline\addlinespace
\addlinespace
\addlinespace
8 & $\left\langle 119,7021,178878,\dots,50400,2520\right\rangle $\tabularnewline\addlinespace
\addlinespace
\addlinespace
9 & $\left\langle 246,30135,1409590,\dots,544320,20160\right\rangle $\tabularnewline\addlinespace
\addlinespace
\addlinespace
$n$ & $\left\langle 2^{n-1}-n-1,\binom{f_{1}}{2},\dots,\frac{n!(n-3)}{4},\frac{\left(n-1\right)!}{2}\right\rangle $\tabularnewline\addlinespace
\bottomrule
\addlinespace
\end{tabular}
\par\end{centering}

\caption{\label{tab:Computational-results-for}Computational results for the
Kalmanson complex.}
\end{table}
For small $n$, the $f$-vector may be computed directly. Results
for $n=4,\dots,9$ are presented in Table \ref{tab:Computational-results-for}.
We now theoretically explain some of these numbers. First we restate
some additional definitions and results from \cite{Bandelt:1992p3044,Dress:1996p3275}
which will prove useful in enumerating the faces of $\mathcal{K}_{n}$.
For the remainder of the section, $S_{i}\in\mathcal{S}(X)$ represents
a split and the identity $S_{i}=\left\{ A_{i},B_{i}\right\} $ is
implicit.
\begin{defn}
\label{def:weakly-compatible}A split system $\mathcal{S}$ is called
\emph{weakly compatible} if for all triples $S_{1},S_{2},S_{3}\in\mathcal{S}$
there do not exist points $a,a_{1},a_{2},a_{3}\in X$ such that $a\in A_{1}\cap A_{2}\cap A_{3}$
and $a_{i}\in A_{j}\iff i=j$.
\end{defn}
\begin{figure}
\noindent \begin{centering}
\includegraphics[scale=1.2]{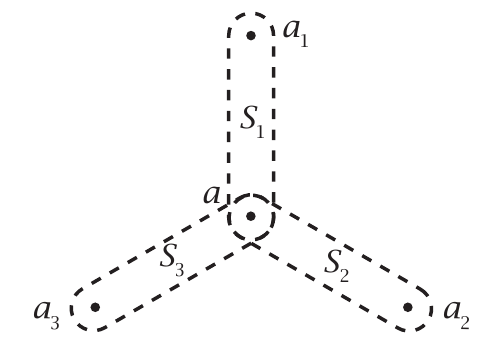}
\par\end{centering}

\caption{\label{fig:weakly-incompatible}A system of splits which is not weakly
compatible.}
\end{figure}
Weak compatibility enforces a sort of convexity condition on $\mathcal{S}$
by requiring that, for any triple of points in the split system, there
is no point which mutually separates them (Figure \ref{fig:weakly-incompatible}).

For any two splits $S_{1},S_{2}$ we define a binary operation $\sqcup$
by $\left\{ A_{1},B_{1}\right\} \sqcup\left\{ A_{2},B_{2}\right\} =\left\{ A_{1}\cap A_{2},B_{1}\cup B_{2}\right\} $.
\begin{lem}
\label{lem:weak_compat}The splits $S_{1}$, $S_{2}$ and $S_{1}\sqcup S_{2}$
are weakly compatible.
\begin{proof}
Let $S_{3}=S_{1}\sqcup S_{2}$. If there exist $a,a_{1},a_{2},a_{3}$
as in Definition \ref{def:weakly-compatible}, then $a_{3}\in A_{3}-(A_{1}\cup A_{2})\neq\emptyset$.
But $A_{3}=A_{1}\cap A_{2}$, a contradiction.
\end{proof}
\end{lem}
\begin{thm}[\cite{Dress:1996p3275}]
\label{thm:cyclic_condition}Let $\mathcal{S}$ be a split system
and let $\mathcal{S}'$ be the split system\[
\mathcal{S}':=\mathcal{S}\cup\Big\{ S_{1}\sqcup S_{2}:S_{1},S_{2}\in\mathcal{S}\text{ \emph{and} }A_{i}\cap B_{j}\neq0,\ i,j\in\left\{ 1,2\right\} \Big\}\]
Then $\mathcal{S}$ is contained in a circular split system if and
only if $\mathcal{S}'$ is weakly compatible.\end{thm}
\begin{cor}
A circular split system is weakly compatible.
\end{cor}

\subsection{Low (Co-)Dimensional Faces}

Enumerating the vertices, edges, ridges and facets of $\mathcal{K}_{n}$
is now straightforward.
\begin{thm}
\label{thm:f-vector}Let $f=(f_{0},\dots,f_{d-1})$ denote the $f$-vector
of $\mathcal{K}_{n}$. Then \begin{align}
f_{0} & =2^{n-1}-n-1\\
f_{1} & =\binom{f_{0}}{2}\label{eq:f_1}\\
f_{d-2} & =\left[\binom{n}{2}-n\right]\times f_{d-1}\label{eq:f_d-2}\\
f_{d-1} & =\frac{\left(n-1\right)!}{2}\label{eq:f_d-1}\end{align}

\begin{proof}
$f_{0}$ counts the number of non-trivial $X$-splits. There are \[
\sum_{k=2}^{n-2}\binom{n}{k}=2^{n}-2n-2\]
binary words on $n$ letters which contain at least two zeros and
two ones. Since each word and its complement correspond to the same
split, we divide by two to obtain $f_{0}$.

Equation \eqref{eq:f_1} asserts that every pair of splits $\mathcal{S}_{1},\mathcal{S}_{2}$
is contained in a circular split system. By Lemma \ref{lem:weak_compat},
the splits $S_{1}$, $S_{2}$ and $S_{1}\sqcup S_{2}$ are weakly
compatible. Then by Theorem \ref{thm:cyclic_condition}, $\left\{ S_{1},S_{2}\right\} $
is a circular split system.

Each facet of $\mathcal{K}_{n}$ corresponds to a circular ordering;
that is, an edge labeling of the regular $n$-gon. Such labelings
are unique up to dihedral symmetry. There are $(n-1)!$ labelings
up to rotation, and half that number when accounting for reflection.
This yields \eqref{eq:f_d-1}.

To prove \eqref{eq:f_d-2}, let $F\subset\mathcal{K}_{n}$ be a facet
spanned by vertices $v_{1},\dots,v_{d}\in X$; without loss of generality
assume the circular ordering corresponding to $F$ is $(1,2,\dots,n)$.
Let $u$ be another vertex distinct from the $v_{i}$, with corresponding
split \[
S_{u}=\Bigg\{\left\{ 1,u_{2},\dots,u_{j}\right\} ,\left\{ u_{j+1},\dots,u_{n}\right\} \Bigg\}\]
Finally, let $i=\min\left\{ i:u_{i}\neq i\right\} $, which exists
by the assumption that $S_{u}$ is not circular with respect to the
given ordering. Now, the splits $S_{u}$ and\begin{align*}
S_{1} & =\Big\{\left\{ u_{i}-1,u_{i}\right\} ,X-\left\{ u_{i}-1,u_{i}\right\} \Big\}\in F\\
S_{2} & =\Big\{\left\{ u_{i},u_{i}+1\right\} ,X-\left\{ u_{i},u_{i}+1\right\} \Big\}\in F\end{align*}
are weakly incompatible: denoting the first blocks of each by $A_{u},A_{1},A_{2}$
we have\begin{align*}
\left\{ u_{i}\right\}  & =A_{u}\cap A_{1}\cap A_{2}\\
1 & \in A_{u}-(A_{1}\cup A_{2})\\
u_{i}-1 & \in A_{1}-(A_{u}\cup A_{2})\\
u_{i}+1 & \in A_{2}-(A_{u}\cup A_{1})\end{align*}
Hence, by contradiction any collection of $d-1$ vertices of $F$
spans a unique face of codimension two. As described in Section \ref{sec:Geometry-of},
each facet contains $\binom{n}{2}-n$ vertices (one for each diagonal
of the $n$-gon.) 
\end{proof}
\end{thm}

\subsection{Triangles}

The computations in Theorem \ref{thm:f-vector} were aided by the
fact that $\mathcal{K}_{n}$ is connected in dimension one and totally
disconnected in codimension one. Enumerating the faces in the remaining
cases is more challenging. To illustrate the issues involved, we demonstrate
how to compute $f_{2}$, the number of triangles in $\mathcal{K}_{n}$. 
\begin{example}
The split system \[
\biggl\{\Big\{\left\{ 1,2\right\} ,\left\{ 3,4,5\right\} \Big\},\Big\{\left\{ 1,3\right\} ,\left\{ 2,4,5\right\} \Big\},\Big\{\left\{ 1,4\right\} ,\left\{ 2,3,5\right\} \Big\}\biggr\}\]
is not weakly compatible, so it is not a triangle of $\mathcal{K}_{n}$.
By contrast, the split system\[
\biggl\{\Big\{\left\{ 1,2\right\} ,\left\{ 3,4,5\right\} \Big\},\Big\{\left\{ 2,3\right\} ,\left\{ 1,4,5\right\} \Big\},\Big\{\left\{ 4,5\right\} ,\left\{ 1,2,3\right\} \Big\}\biggr\}\]
is circular with respect to two orderings: $\left(1,2,3,4,5\right)$
and $(1,2,3,5,4)$. It is therefore a triangle of \emph{$\mathcal{K}_{n}$}
which is contained in two facets.
\end{example}
Our main tool for computing $f_{2}$ will be Corollary \ref{cor:circular iff c1r},
in conjunction with a structure theorem of \cite{Tucker:1972p5021}
which completely characterizes C1R matrices.
\begin{defn}
Let $M$ be a matrix. The \emph{configuration of }$M$ is the set
of matrices obtained by permuting the rows and/or columns of $M$
(not necessarily by the same permutation.)\end{defn}
\begin{example}
The configuration of the $2\times2$ identity matrix is the set \[
\left\{ \begin{pmatrix}1 & 0\\
0 & 1\end{pmatrix},\begin{pmatrix}0 & 1\\
1 & 0\end{pmatrix}\right\} \]
\end{example}
\begin{thm}[\cite{Tucker:1972p5021}]
\label{thm:forbidden}A binary matrix $M$ is C1R if and only if
it does not contain as a submatrix any configuration of $M_{\text{I}_{n}},M_{\text{II}_{n}},M_{\text{III}_{n}},M_{\text{IV}},M_{\text{V}}$,
$1\leq n<\infty$, where \begin{flalign*}
M_{\text{I}_n} &= \kbordermatrix{ 	& c_1 & c_2 & c_3 & \cdots & c_n & c_{n+1} & c_{n+2} \cr
				r_1	& 1   &  1  &  0  & \cdots &  0  &  0      &   0     \cr
				r_2 & 0   &  1  &  1  & \cdots &  0  &  0      &   0     \cr
				\vdots & \vdots & \vdots & \vdots & \ddots & \vdots & \vdots & \vdots \cr
				r_n & 0 & 0 & 0 & \cdots & 1 & 1 & 0 \cr
				r_{n+1} & 0 & 0 & 0 & \cdots & 0 & 1 & 1 \cr
				r_{n+2} & 1 & 0 & 0 & \cdots & 0 & 0 & 1 
}& M_{\text{IV}} &= \kbordermatrix{ 	& c_1 & c_2 & c_3 & c_4 & c_5 & c_6 \\
				r_1 & 1 & 1 & 0 & 0 & 0 & 0\\
				r_2 & 0 & 0 & 1 & 1 & 0 & 0\\
				r_3 & 0 & 0 & 0 & 0 & 1 & 1\\
				r_4 & 0 & 1 & 0 & 1 & 0 & 1
}& \\
M_{\text{II}_n} &= \kbordermatrix{ 	& c_1 & c_2 & c_3 & \cdots & c_n & c_{n+1} & c_{n+2} & c_{n+3} \cr
				r_1	& 1   &  1  &  0  & \cdots &  0  &  0      &   0     & 0 \cr
				r_2 & 0   &  1  &  1  & \cdots &  0  &  0      &   0     & 0 \cr
				\vdots & \vdots & \vdots & \vdots & \ddots & \vdots & \vdots & \vdots & \vdots \cr
				r_n & 0 & 0 & 0 & \cdots & 1 & 1 & 0 & 0 \cr
				r_{n+1} & 0 & 0 & 0 & \cdots & 0 & 1 & 1 & 0 \cr
				r_{n+2} & 1 & 1 & 1 & \cdots & 1 & 1 & 0 & 1 \cr
				r_{n+3} & 0 & 1 & 1 & \cdots & 1 & 1 & 1 & 1
}& M_{\text{V}} &= \kbordermatrix{ 	& c_1 & c_2 & c_3 & c_4 & c_5 \\
				r_1 & 1 & 1 & 0 & 0 & 0 \\
				r_2 & 1 & 1 & 1 & 1 & 0 \\
				r_3 & 0 & 0 & 1 & 1 & 0 \\
				r_4 & 1 & 0 & 0 & 1 & 1 
}\displaybreak \\
M_{\text{III}_n} &= \kbordermatrix{ 	& c_1 & c_2 & c_3 & \cdots & c_n & c_{n+1} & c_{n+2} & c_{n+3} \cr
				r_1	& 1   &  1  &  0  & \cdots &  0  &  0      &   0     & 0 \cr
				r_2 & 0   &  1  &  1  & \cdots &  0  &  0      &   0     & 0 \cr
				\vdots & \vdots & \vdots & \vdots & \ddots & \vdots & \vdots & \vdots & \vdots \cr
				r_n & 0 & 0 & 0 & \cdots & 1 & 1 & 0 & 0 \cr
				r_{n+1} & 0 & 0 & 0 & \cdots & 0 & 1 & 1 & 0 \cr
				r_{n+2} & 0 & 1 & 1 & \cdots & 1 & 1 & 0 & 1
} 
\end{flalign*}

\end{thm}
In the specific case of $f_{2}$, where we are counting $3\times n$
matrices, only two forbidden submatrices pertain:\[
M_{I_{1}}=\begin{pmatrix}\begin{array}{rrr}
1 & 1 & 0\\
0 & 1 & 1\\
1 & 0 & 1\end{array}\end{pmatrix}\qquad\text{and}\qquad M_{III_{1}}=\begin{pmatrix}\begin{array}{rrrr}
1 & 1 & 0 & 0\\
0 & 1 & 1 & 0\\
0 & 1 & 0 & 1\end{array}\end{pmatrix}\]

\global\long\def\col{\operatorname{col}}
For a matrix $M$, we write $\col\left(M\right)$ to denote the set
of column vectors of $M$. Let $I=\col\left(M_{I_{1}}\right)$ and
$III=\col\left(M_{III_{1}}\right)$. Note that $I$ and $III$ are
{}``closed'' under the operation of row permutation. Hence, by Corollary
\ref{cor:circular iff c1r} and Theorem \ref{thm:forbidden},\[
[M]\in\mathcal{Q}_{3}\iff\left|\col(M)\cap I\right|<3\text{ and }\left|\col(M)\cap III\right|<4\]

Accordingly, let \begin{equation}
F_{i,j}=\left\{ \left[M\right]\in\mathcal{Q}_{3}:\left|\col(M)\cap I\right|=i\text{ and }\left|\col(M)\cap III\right|=j\right\} \label{eq:Fij}\end{equation}
Then \begin{equation}
f_{2}=\left|\mathcal{Q}_{3}\right|=\sum_{\substack{0\leq i\le2\\
0\leq j\leq3}
}\left|F_{i,j}\right|\label{eq:sum Fij}\end{equation}

Enumerating $F_{i,j}$ involves carefully counting the number of classes
of $\mathcal{Q}_{3}$ while keeping track of how many columns from
the sets \begin{align*}
I & =\left\{ \begin{pmatrix}1\\
1\\
0\end{pmatrix},\begin{pmatrix}1\\
0\\
1\end{pmatrix},\begin{pmatrix}0\\
1\\
1\end{pmatrix}\right\}  & III & =\left\{ \begin{pmatrix}1\\
0\\
0\end{pmatrix},\begin{pmatrix}0\\
1\\
0\end{pmatrix},\begin{pmatrix}0\\
0\\
1\end{pmatrix},\begin{pmatrix}1\\
1\\
1\end{pmatrix}\right\} \end{align*}
appear in each equivalence class.

\subsubsection{\label{sub:Sample-Calculation:}Sample Calculation: $\left|F_{0,3}\right|$}

The counting argument is straightforward but tedious. We illustrate
the calculation of $\left|F_{0,3}\right|$; the remaining cases are
similar and are proven in \cite{terhorst}.

Let $\mathbb{P}_{n}$ denote the set of ordered partitions of the
integer $n$. That is, for a $k$-tuple $(x_{1},\dots,x_{k})$ we
have\[
(x_{1},\dots,x_{k})\in\mathbb{P}_{n}\iff\sum_{i=1}^{k}x_{i}=n\text{ and }x_{i}\geq1\text{ for all }i\]

To simplify the notation we take summation over $\mathbb{P}_{n-1}$
for granted wherever there is no chance of confusion: instead of e.g.\[
\sum_{\substack{(a,b,c,d)\in\mathcal{P}_{n-1}\\
a>1}
}\binom{n-1}{a,b,c,d}\]
we will simply write\[
\sum_{\substack{a>1}
}\binom{n-1}{a,b,c,d}\]

Now let $[M]\in F_{0,3}$. We consider two cases.
\begin{enumerate}
\item First, if $\left(1,1,1\right)^{T}\notin\col\left(M\right)\cap III$
then $M$ is of the form \[
\begin{tikzpicture}[mymatrixenv]     
\matrix[mymatrix] (m) 
{ 0 & 0 & 0 & 1 & 1 & 0 & 0 & 0\\
  0 & 0 & 0 & 0 & 0 & 1 & 1 & 1\\
  0 & 1 & 1 & 0 & 0 & 0 & 0 & 0\\
};
\mymatrixbracetop{2}{3}{$a$}     
\mymatrixbracetop{4}{5}{$b$}     
\mymatrixbracetop{6}{8}{$c$}
\end{tikzpicture} 
\]where $a,b,c$ count the instances of the columns in $III-\left\{ \left(1,1,1\right)\right\} $.
To prevent the occurrence of a trivial split (row containing $<2$
ones) we require $\min\left(a,b,c\right)>1$. Hence there are \begin{equation}
(1/6)\hspace{-0.2in}\sum_{\min(a,b,c)>1}\negthickspace\binom{n-1}{a,b,c}\label{eq:f0300}\end{equation}
such classes, where the factor of $1/6$ reflects the fact that each
arrangement is equivalent to six others obtained by permuting the
labelings $a,b,c$. We see that this counts the number ways of arranging
$n-1$ ones into three unlabeled rows where each must contain at least
two ones.%
\footnote{This number is also known as an associated Stirling number of the
second kind, cf. A000478 \cite{oeis}.%
}\\
\\
We also have the possibility that $M$ contains additional zero
columns:\[ 
\begin{tikzpicture}[mymatrixenv]     
\matrix[mymatrix] (m) 
{ 0 & 0  & 1 & 1 & 0 & 0  & 0 & 0\\
  0 & 0  & 0 & 0 & 1 & 1  & 0 & 0\\
  0 & 0  & 0 & 0 & 0 & 0  & 1 & 1\\
};
\mymatrixnakedtop{2}{2}{$d$}     
\mymatrixbracetop{3}{4}{$a$}     
\mymatrixbracetop{5}{6}{$b$}
\mymatrixbracetop{7}{8}{$c$}
\end{tikzpicture} 
\]By an entirely analogous argument we count\begin{equation}
(1/6)\hspace{-0.2in}\sum_{\min(a,b,c)>1}\negthickspace\binom{n-1}{a,b,c,d}\label{eq:f0301}\end{equation}
such classes. 
\item In the second case we have $\left(1,1,1\right)^{T}\in\col\left(M\right)\cap III$.
These are matrices of the form\[ 
\begin{tikzpicture}[mymatrixenv]     
\matrix[mymatrix] (m) 
{ 0 & 1 & 1 & 1 & 0 & 0 & 1 & 1\\
  0 & 1 & 1 & 1 & 0 & 0 & 0 & 0\\
  0 & 1 & 1 & 1 & 1 & 1 & 0 & 0\\ 
};
\mymatrixbracetop{2}{4}{$c$}     
\mymatrixbracetop{5}{6}{$a$}
\mymatrixbracetop{7}{8}{$b$}
\end{tikzpicture} 
\]We must have that the number of $(1,1,1)^{T}$ columns is greater
than one, or else we have a trivial split. Also, each arrangement
is equivalent to the one obtained by swapping the columns labeled
$a$ and $b$ and permuting their respective rows. We therefore count
\begin{equation}
(1/2)\sum_{c>1}\binom{n-1}{a,b,c}\label{eq:f0210}\end{equation}
such classes. We also obtain \begin{equation}
(1/2)\sum_{c>1}\binom{n-1}{a,b,c,d}\label{eq:f0211}\end{equation}
classes by allowing for the presence of additional zero columns.
\end{enumerate}
Now let $S_{n,k}$ denote the Stirling number of the second kind.
We will make use of the identity \begin{equation}
\sum_{(x_{1},\dots,x_{k})\in\mathbb{P}_{n}}\binom{n}{x_{1},\dots,x_{k}}=k!\cdot S_{n,k}=:M(n,k),\label{eq:numsurjections}\end{equation}
to obtain simple formulas for equations \eqref{eq:f0300}--\eqref{eq:f0211}.
(One interpretation of $M(n,k)$ is that it counts the number of surjections
from an $n$-set onto a $k$-set, see e.g.\emph{\ }\cite[Ch. 1]{Aigner:2007fk}.)
By \eqref{eq:numsurjections}, symmetry and the inclusion-exclusion
principle we have \begin{multline}
\sum_{\min(a,b,c)>1}\negthickspace\binom{n-1}{a,b,c}=\\
M(n-1,3)-3\sum_{a=1}\negthickspace\binom{n-1}{a,b,c}+3\negthickspace\sum_{a=b=1}\binom{n-1}{a,b,c}-\negthickspace\sum_{a=b=c=1}\negthickspace\binom{n-1}{a,b,c}\label{eq:simplify}\end{multline}

The second term of \eqref{eq:simplify} counts the number of surjections
from a set of cardinality $n-1$ onto the set $\left\{ a,b,c\right\} $
such that a unique element maps to $a$. There are $\left(n-1\right)\times M(n-2,2)$
such maps. Similarly, the third and fourth terms represent $2!\times\binom{n-1}{2}\times M(n-3,1)=(n-1)(n-2)$
and $3!\times(n-1,3)\times M(n-4,0)=0$ maps respectively. We therefore
conclude \[
\sum_{\min(a,b,c)>1}\negthickspace\binom{n-1}{a,b,c}=M(n-1,3)-3(n-1)M(n-2,2)+3(n-1)(n-2)\]
By the same arguments, \begin{multline*}
\sum_{\min(a,b,c)>1}\negthickspace\binom{n-1}{a,b,c,d}=M(n-1,4)-3(n-1)M(n-2,3)+\\
3(n-1)(n-2)M(n-3,2)-(n-1)(n-2)(n-3)\end{multline*}

For \eqref{eq:f0210} we note that \[
\sum_{c>1}\binom{n-1}{a,b,c}=M(n-1,3)-\sum_{c=1}\binom{n-1}{a,b,c}=M(n-1,3)-(n-1)\times M(n-2,2)\]
Similarly, we can rewrite \eqref{eq:f0211} as $M(n-1,4)-(n-1)\times M(n-2,3)$.

\subsubsection{General Formula}

Repeating these counting arguments for the remaining $\left|F_{i,j}\right|$
yields the following formula.
\begin{thm}
\label{thm:triangles}Let $t=n-1$. The number of triangles in $\mathcal{K}_{n}$
is\begin{multline*}
(1/6)(t-2)(t-1)t+2(t-1)t\left[1+M(t-2,2)\right]-\\
5t\cdot M(t-1,2)-8t\cdot M(t-1,3)-2t\cdot M(t-1,4)+\\
\left(19/6\right)M(t,3)+\left(55/6\right)M(t,4)+7M(t,5)+2M(t,6)\end{multline*}

\begin{proof}
See \cite{terhorst}.
\end{proof}
\end{thm}
The first ten entries of this sequence, $n=4,\dots,13$ are \[
0,90,1755,19950,178878,1409590,10270585,71110930,475443364,3100707610,\dots\]
We have verified this formula computationally up to $n=10$ (the largest
$n$ for which the calculations terminated) using the mathematics
software \textsf{SAGE} \cite{sage}. Source code for this and related
$f$-vector calculations may be downloaded from: \url{https://github.com/terhorst/kalmanson}.

\section{\label{sec:Conclusion}Conclusion}

In this paper we analyzed the combinatorics of the Kalmanson complex.
We show how this complex arises in split theory, optimization and
phylogenetics. We gave a simplified proof of the equivalence of Kalmanson
and circular decomposable metrics based on polyhedral geometry and
our interpretation of $\mathcal{K}_{n}$ as a simplicial complex of
splits.

Subsequently, our main focus was to enumerate its faces. We demonstrated
that the complex is totally connected along edges, and totally disconnected
along ridges. We then gave a formula for enumerating its triangles,
using a forbidden substructure characterization along with some basic
counting principles.

At present we do not have a way to generalize this method to faces
of arbitrary dimension. The next case of tetrahedra ($k=4$) becomes
considerably more difficult, as there are now 7 avoided Tucker matrices
to consider: $M_{I_{1}},M_{I_{2}},M_{II_{1}},M_{III_{1}},M_{III_{2}},M_{IV},M_{V}$.
The connection to the Tucker theorem suggests a possible application
of results on avoided configurations (see \cite{anstee} for a survey),
but most results in that literature are of an extremal, as opposed
to enumerative, variety. In \cite{Spinrad:2003fk} some matrices avoiding
small configurations are counted, but we are not aware of a general
method of enumerating matrices which avoid configurations of arbitrary
dimensions. We view this as in interesting problem in enumerative
combinatorics which merits further study.

\newpage{}

\bibliographystyle{amsplain}
\bibliography{../kalmanson_papers,../other_papers}

\end{document}